\documentstyle[12pt]{article}\topmargin=-2cm\textheight=220mm
\begin{document}
\begin{center}{\large\bf New Geometric Framework for $SU(2)$ Gauge
Theory}\\
Zafar Ya.\,Turakulov\\ {\it Institute of Nuclear Physics, Ulugbek,
Tashkent, 702132 Uzbekistan (e-mail: zafar@suninp.tashkent.su)}\end{center}
\begin{abstract}An explicit model of fiber bundle with local fibers
being disinct copies of vector 3-space is introduced. They are endowed
with frames which are used as local isotopic ones. The field local of
isotopic frames is considered as gauge field itself while the form of
gauge connections is derived from it. A covariant equation for the field
of local frames is found. It is shown that Yang-Mills equation follows
from it, but variety of solutions of the new equation is highly reduced in
such that no ambiguities (Yang-Wu and vacuum ones) arise. It is shown that
Lagrangian for the field gives non-zero trace for the stress-energy tensor
and zero value for spin of the field of plane wave. Some new solutions for
the fields of punctual source and spherical wave are found.\end{abstract}
\newpage\section{Introduction}

Fiber bundles are known to provide adequate geometric models for
classical gauge theories in which classical gauge fields are described
by connections on the bundles whose sections are scalar fields with
isospin (Wu and Yang, 1975; Trautman, 1979). Ordinarily, one introduces
such a section as a column of three scalar fields and calls them isotopic
components of the field (Actor, 1979, Ryder, 1985). Isotopic rotations
are well-defined on the space of columns and this definition allows
introducing the notion of gauge field as small rotations corresponding to
small translations in the space-time. Though this approach is generally
accepted it leads to the following difficulty.

First, since local isotopic spaces are vector ones they are to be endowed
with a frame for each of them. Second, connection, which, by definition,
connects the frames, can be derived in the same way as it is done in
Riemannian geometry. Third, since some three fields are called isotopic
components they are to be projections of an isovector onto elements of the
frame. And, last, the connection must give in general non-zero curvature.
In fact, these four requirements are inconsistent.

Indeed, if the frames are introduced as columns they belong to one and the
same space $R^3$, not to different copies of a vector space. And since
they are elements of one and the same space the bundle obtained this way
has by construction, identically zero curvature. Triplets of such columns
forming local frames can be parallelized by some appropriate field of
local isotopic rotations which gauges out the corresponding connection. The
goal of this work is to construct an explicit model of fiber bundle such
whose local fibers are different copies of a vector space each of which can
be endowed with a frame, the connection can be obtained from these frames
and gives non-zero curvature.

It is expected that if this is made the theory suffers some changes.
Indeed, if a field of local isotopic frames is introduced and the
connection can be derived from the local frames and, hence, it is not
longer an independent variable. Since it is the field of local frames which
predetermines explicitly all the rest, i.\,e.\,connection and curvature,
and thus specifies the gauge field itself. Since it is not evident that
some bundle of frames exists for any given 1-form accepted as connection,
variety of connections derived from local frames may be highly reduced.
As well-known, in Riemannian geometry connections specified arbitrarily may
contain torsion while those found from local frames constitute a special
class of torsion-free connections. If something similar takes place in
other fiber bundles i.\,e.\,their connections contain in general two terms
of different natureas connection and torsion are, there exists a problem to
separate them and introduce a special field equation for each of them.
Otherwise, either only one of them is to be considered in the theory or
the extra one obeying no any equation will destroy univalence of field
produced by a certain source.

If the field is represented by local frames and connection is derived
from them there is no extra fields, but the field equation is now to be
written for the frames and it obviously differs from the standard the
Yang-Mills equation. Such an approach looks to fit C.\,N.\,Yang and
T.\,T.\,Wu's claim that it is firber bundle which represents the field
itself (Wu and Yang, 1975; see also Socolowski, 1991). However the
generally accepted point of view reads that one needs fiber bundles
`when studying topologically non-trivial objects' (Trautman, 1979). In
fact, one deals with fiber bundle whenever there is a gauge field
regardless of its topological properties.

If, like all field equations of this sort, the equation for local isotopic
frames is of the second order, it contains the curvature components not
their derivatives. Thus, it is expected to be somewhat similar to the
Einstein equation. As for the field properties like mass and spin, they
are to be studied specially.\section{Natural models of fiber bundles}

As was pointed out above it is not sufficient to write down a column of
three numbers and call them isotopic vector in given spacetime point. To
introduce genuine local frames for a fiber bundle on needs collection of
spaces forming the bundle in which the frames can be introduced. Thus, the
first problem is to find an explicit representation of a bundle in which
all the fibers are really different copies of a space. It seems that the
only way to obtain such a construction is to find among known ones.

Examples of bundles with fibers given as different copies of a vector
space are well-known. One of them is co-tangent bundle $\Lambda^1M$ over
a Riemannian manifold $M$ whose sections are ordinary fields of 1-forms
on $M$ (Warner, 1983). If $\{x^i\}$ is a coordinate system on $M$ each
fiber is endowed with a frame $\{dx^i\}$ which can be normailzed by
appropriate local linear transformations. Geometry of the bundle is
predetermined by that of its base $M$ and its frames are consistent with
the only non-trivial connection found from them by solving the first
structure equation as a system of algebraic equations for theconnection
components. So are the bundles $\Lambda^pM$ of $p-$forms on the manifold.
All of them have one and the same structure group which is Lorentz group
(here and below $M$ is (3+1)-dimensional).

In order to introduce bundles with other structure groups one can consider
so-called subbundles of $\Lambda^pM$, i.\,e.\,bundles whose typical fibers
are subspaces of local fibers $\Lambda^p_x,x\in M$. For example a certain
choice of 3-dimensional spacelike hyperplanes in local co-tangent spaces
$\Lambda^1_x$ specifies a new fiber bundle with $SO(3)$ as the structure
group. Another model with this group occures when chosing 3-dimensional
subspaces in local spaces $\Lambda^2_x$. An interesting example of a
bundle with $SU(3)$ as the structure group can, in principle, be composed
by taking complex combinations of 2-forms in $\Lambda^2_x$ and their duals
and introducing Hermitean metric for the spaces of such combinations
requiring the Lorentz transformations to enter as a $SU(2)$ sugroup of the
structure group. Perhaps, no other structure groups but $U(1)$ can be
obtained this way. In this work we restrict ourselves with considering one
case with local fibers being 3-dimensional subspaces in local spaces
$\Lambda^2_x$ and $SO(3)$ as the structure group.

Since local fibers constructed this way are subspaces of $\Lambda^2_x$'s
the corresponding structure group is a subgroup of the Lorentz group.
However, we seek to employ the subbundle to describe gauge fields and,
hence, rotations of local fibers as local gauge transformation that are
to be inaffected by any coordinate transformations in the spacetime. It is
possible because it suffices to have only the spaces over spacetime points
as local fibers. It will be shown below how they may be used as models of
local isotopic spaces.\section{Local frames, gauge transformations and
connection}

Consider the fiber bundle $\Lambda^2M$ of 2-forms over the spacetime $M$
and select a 3-dimensional subspace $I_x$ in each local fiber
$\Lambda^2_x$. As was pointed out above we use only the collection of
spaces $I_x$ that constitute a vector bundle $IM$. As subspaces of
$\Lambda^2_x$' they have a certain metric inherited from these spaces.
However, it is unnecessary to accept all their properties like this
metric. They themselves form the structure of fiber bundle as a whole,
and metric on each of them can be introduced arbitrarily.

This may be done as follows. Let $\{\pi^a\}_x,\enskip\pi^a\in I_x$ be a
triplet of 2-forms. A new Euclidean metric for spaces $I_x$ can be
introduced simply by postulating that $\pi^a$'s constitute orthonormalized
frames in the spaces. Thus, in general, $\pi^a$'s are neither orthogonal
nor nomalized with respect to the natural metric of $\Lambda^2_x$'s but
they are orthonormalized due to the new metric introduced thereby.
Hereafter the spaces $I^a_x$ and frames $\pi^a$'s are considered as models
of local isotopic spaces and as local isotopic frames respectively. The
isotopic vectors $\pi^a$ are invariant under Lorentz transformations
(though their coordinate components obey the usual transformation law). In
other words, coordinate transformations in $M$ change their Lorentzian
components, but the vectors themselves remain the same. Since isotopic
vectors are referred to the local frames $\{\pi^a\}_x$ and, hence, they
have no coordinate components, their components are inaffected by any
coordinate transformations in the spacetime. Thus, we have composed a model
of fiber bundle whose fibers are different copies of a vector space by
construction.

Now, local gauge transformations can be introduced as linear
transformations conserving orthogonality and normalization of the isotopic
frames $\{\pi^a\}_x$. They are defined as follows. Let $S_x$ be a linear
operator that changes the metric of space $I_x$ into the natural metric
induced on it from $\Lambda^2M$ and $R$ be a usual rotation in the space
$\Lambda^2_x$. Then, linear operators of the form $S_x^{-1}R_xS_x$ local
isotopic frames $\{\pi^a_x\}$ orthonormalized. At the same time they
constitute the group $SO(3)$ which, therefore, can be considered as a model
of the group of gauge transformations.

Infinitesimal parameters of a small gauge transformation $dg^a$
referred to the chosen local isotrpic frame can now be introduced by the
standard procedure: if $\{\pi^a\}_x$ and $\{\pi'^a\}_x$ are two frames
connected by a small gauge trnsformation then
$$\pi'^a_x=\pi^a_x+\varepsilon^a{}_{bc}dg^b\pi^c.$$This definition allows
introducing a connection that really connects isotopic frames established
in two neighbouring spacetime points via the first structure equation
\begin{equation}
d\pi^a+\varepsilon^a{}_{bc}\alpha^b\pi^c=0\label{str1}\end{equation}
where the connection appeared as 1-form $\alpha^b$. Hereafter this
equation is used as the definition of the connection. The definition of
the connection form just proposed is univalent. Indeed, as was shown by
R.\,Roskies (Roskies, 1977), M.\,Calvo (Calvo, 1977) and M.\,Halpern
(Halpern, 1977) the Yang-Wu ambiguity of the vector potential does not
arise if the algebraic system for components of $\alpha$ is well-defined,
i.\,e.\,if it consists of 12 equations. As the gauge connection defined all
the rest can be built in the ordinary manner (see Actor, 1979 or Ryder,
1985 for details).\section{The field equations and Lagrange formalism}

Now compose the field equation for the field of local frames $\pi^a$
assuming it to be covariant and forminvariant with respect to gauge
transformations as we introduced them i.\,e.\,ordinary ortogonal
transformations of the frame. Therefore, the equation is to be
composed of the frame itself and the curvature 2-form $K^a$ which is
\begin{equation}K^a=d\alpha^a+\frac{1}{2}\varepsilon^a{}_{bc}
\alpha^b\wedge\alpha^c\label{curvature}\end{equation}The only
manifestly covariant equation of the second order on $\pi^a$'s composed
of these two elements is\begin{equation}K^a-m^2a\pi^a=0\label{fieldeq}
\end{equation}where $m$ is a constant. Apparently, such a field
equation is consistent with Bianchi identities because $\pi^a$'s
satisfy the first structure equation (\ref{str1}). Denoting, as usual,
the operator acting on $\pi^a$ in the equation $D\pi^a$ and calling it
covariant exterior derivative we obtain the standard Yang-Mills
equation in the form\begin{equation}D{}^\ast\pi^a=I^a\label{yme}
\end{equation} with $I^a$ being the 3-form of conserving current
$${}^\ast I^a\equiv J_i^adx^i,\enskip DI=0.$$ It should be noted that
normalization of the current as it is defined here differs from the
standard one by the constant factor $m^2$.. A similar approach to
Yang-Mills equations has been proposed by L.\,Castillejo and M.\,Kugler
in their work (Castillejo and Kugler, 1980).

Now we pass temporarily to tensorial denotions to describe Langrangian
formalism for the equations just proposed. All the equations can be
derived in standard manner from the following Lagrangian density
\begin{equation}L=(\partial_iA_j+
\frac{1}{2}\varepsilon^a{}_{bc}A^b_iA^c_j)\pi^{ij}_a-\frac{m^2}{2}
\pi_{ij}^a\pi_{ij}^a+J^i_aA_i^a\label{lagrangian}\end{equation}where
$\alpha^a\equiv A_i^adx^i$, $\pi^a\equiv\pi^a_{ij}dx^i\wedge dx^j$ and
$J^i_a$ stands for currents of an external source. Indeed, since the
Lagrangian does not include derivatives of $\pi^a$'s one of the
Euler-Lagrange equations reads simply $$\frac{\partial
L}{\partial\pi^a_{ij}}=0$$ that is exactly the equations
(\ref{curvature}) and (\ref{fieldeq}). Now the equation (\ref{str1})
follows from Bianchi identities for $A^a_i$'s. To derive another
Euler-Lagrange equation we evaluate the following derivatives:
$$\frac{\partial L}{\partial A^d_k}=
\frac{1}{2}\varepsilon^a{}_{bc}(\delta^b_d\delta^k_iA^c_j+
\delta^c_d\delta^k_jA^b_i)\pi^{ij}_a+J^k_d=$$
$$\frac{1}{2}(\varepsilon^a{}_{dc}A^c_j\pi^{kj}_a+
\varepsilon^a{}_{bd}A^b_i\pi^{ik}_a)+J^k_d,\quad
\frac{\partial L}{\partial\partial_k A^d_l}=\nabla_k\pi^{kl}_d,$$
where $\nabla_i$ denotes the usual Riemannian covariant derivative.
It is seen that the Euler-Lagrange equation in question coincides with
the Yang-Mills equation (\ref{yme}).

As the explicit form of the Lagrangian is known it is possible to
determine the main characterisitcs of the field i.\,e.\,its mass, spin
and isospin. Spin and isospin of the field will be evaluated for a plane
wave in one of the following sections together with the corresponding
solution and the field mass can be found regardless of the form of
solution directly from that of Lagrangian. In the manifestly covariant
form the field stress-energy tensor is\begin{equation}
T^m_n=\frac{\partial L}{\partial K^a_{mi}}K^a_{ni}-\delta^m_nL=
\pi^a_{mi}K^{ni}_a-\delta^m_nL,\label{set}\end{equation}whereas, as
follows from the equations (\ref{fieldeq}) and (\ref{lagrangian}) the
source-free Lagrangian is equal to $\frac{1}{2}m^2\pi_{ij}^a\pi^{ij}_a$.
Therefore, the trace of the stress-energy tensor is non-zero. It is equal
to $-m^2\pi_{ij}^a\pi^{ij}_a$ that means that the field is massive. On the
other hand, the factor $m$ cannot be put zero because this would erase
mathematical structure of the theory, particularly, the links between
the field of local isotopic frames and the Yang-Mills equation (\ref{yme})
as it appears in the above considerations.

Apart from univalence of links between the source and vector potential
the mathematical framework for theory of $SU(2)$ gauge field proposed
has another important feature. As seen from the equations (\ref{str1})
and (\ref{fieldeq}) the Yang-Mills equation follows from them, thus,
any solution of (\ref{fieldeq}) satisfies the Yang-Mills equation.
In fact, introducing the field of local frames is nothing but order
lowering operation for the Yang-Mills equation. In the following sections
it will be demonstrated that, unlike original Yang-Mills equation, the
system (\ref{str1}-\ref{fieldeq}) is soluble simply by the method of
variables separation.\section{The field of a pointlike non-Abelian charge}

The notion of spherical symmetry has been generalized to the case of
Yang-Mills field by considering the $SU(2)$ group of combined
transformations with $\vec L+\vec T$ as their generators where $\vec L$ is
generator of spatial and $\vec T$ is that of local rotations of isotopic
space (Wu and Wu, 1974; Wilkinson and Goldhaber, 1977). This symmetry has
been used for obtaining an exact solution for a source-free field first by
T.\,T.\,Wu and C.\,N.\,Yang (Wu and Yang, 1975). Spherically-symmetric
solution they found has only spatial components and describes a pointlike
magnetic monopole . It was extended by J.\,P.\,Hsu and E.\,Mac to the case
of a field with non-trivial time-like components which are. Each of three
time components of the vector potential is proportional to one of Cartesian
coordinates in the space, thus, none of them is spherically symmetric (Hsu
and Mac, 1977).

The approach proposed above gives other solutions. To obtain them we
introduce the field of local frames $\{\pi^a\}$ in spherical coordinates
$\{t,r,\theta,\varphi\}$in the form\begin{equation}\pi^1=-P(r)dt\wedge
dr+Q(r)r^2\sin\theta d\theta\wedge d\varphi\label{pisph}\end{equation}
$$\pi^2=p(r)dt\wedge d\theta+q(r)\sin\theta d\varphi\wedge dr$$
$$\pi^3=p(r)\sin\theta dt\wedge d\varphi+q(r)dr\wedge d\theta.$$
and it is seen that the conjugation acts as transmutation
\begin{equation}P\rightarrow -Q,\enskip Q\rightarrow P,\enskip
p\rightarrow -q,\enskip q\rightarrow p.\label{hodge}\end{equation}
To solve the first structure equation (\ref{str1}) we evaluate exterior
derivatives of $\pi^a$'s:
$$d\pi^1=(r^2Q)'\sin\theta dr\wedge d\theta\wedge d\varphi$$
$$d\pi^2=-p'dt\wedge dr\wedge d\theta-
q\cos\theta dr\wedge d\theta\wedge d\varphi$$
$$d\pi^3=p'\sin\theta dt\wedge d\varphi\wedge dr-
p\cos\theta dt\wedge d\theta\wedge d\varphi.$$Inserting this into the
first structure equation (\ref{str1}) and source-free Yang-Mills
equation (\ref{yme}) that can be rewritten as the latter in which the
substitutions (\ref{hodge}) are made, with the connection 1-form
$\alpha^a$:\begin{equation}
\alpha^1=\Phi(r)dt-\cos\theta d\varphi;\quad\alpha^2=A(r)\sin\theta
d\varphi;\quad\alpha^3=A(r)d\theta.\label{alphasph}\end{equation}gives
the following:\begin{equation}(r^2Q)'=2Aq,\quad p'=AP+\Phi q,\quad
(r^2P)'=2Ap,\quad q'=AQ-\Phi p.\label{four}\end{equation}
The curvature 2-form corresponding to the connection (\ref{alphasph})
is $$K^1=-\Phi'r^2dt\wedge dr
-r^{-2}(1-A^2)\sin\theta d\theta\wedge d\varphi$$
$$K^2=A\Phi dt\wedge d\theta-A'\sin\theta d\varphi\wedge dr;$$
$$K^3=A\Phi\sin\theta dt\wedge d\varphi-A'dr\wedge d\theta.$$
Apparently, if the 2-forms $\pi^a$ constitute an orthonormal frame in
the spacetime the curvature is zero. As seen from the form of $\pi^a$
(\ref{pisph}) such a frame corresponds to $P=Q=1,\enskip p=q=r$ and then
$\Phi=0\enskip A=1$. Therefore, if the field is asymptotically trivial,
i.\,e.\, has no strenth, we are to put $A=1$ at infinity. Inserting this
into the field equation (\ref{fieldeq}) gives another quartet of equations:
$$\Phi'=m^2P,\quad m^2Q=-r^{-2}(1-A^2),\quad m^2p=A\Phi,\quad m^2q=A'.$$
Now it is seen that the first two equtions (\ref{four}) are not
independent because they can be derived from the rest ones:
$$(r^2Q)'-2Aq=m^{-2}[(1-A^2)'-AA']\equiv0$$ $$
p'-AP-\Phi Q=m^{-2}[(A\Phi)'-A\Phi'-A'\Phi]\equiv0.$$ These two equations
are satisfied authomatically because they are nothing but Bianchi
identities. Thus, there are only six equations which would be independent.
However, the last quartet of equations specifies only the form of
connection in terms of the coefficients in the definition of the field
of local frame (\ref{pisph}). Thus, it remains to solve the well-known
system of ordinary differential equations (Hsu and Mac, 1977):
\begin{equation}A''-r^{-2}A(1-A^2)-A\Phi^2=0,\end{equation}
$$(r^2\Phi')'-2A^2\Phi=0$$which could also be derived from the Yang-Mills
equation for the connection (\ref{alphasph}). In our case the connection
is found from the first structure equation (\ref{str1}) for localframes
\ref{pisph}) and the Yang-Mills form $D^\ast K$ is annulated by the
independent part of the system of equations (\ref{four}). Due to
univalence of our definition of connection the expressions obtained
represent the unique form of field produced by pointlike non-Abelian
charge. The equations are not soluble in analytical form however,
asymptotical behaviour of the field in question at large $r$ may be
found easily. Indeed, since at infinity the function $A$ is equal to unit
the field behaves asymptotically as $\Phi(r)=C_1r^{-2}$,\enskip$A(r)=1+
C_2r^{-1}$. This solution was obtained in our work (Turakulov 1, 1995).
\section{The field of plane wave}

For the next example we consider the field equations in the standard
coordinate system of circular cylinder $\{t,z,\rho,\varphi\}$. It is
convenient to introduce first a parameter of Lorentz transformations
that allows transferring the solution to be found into an arbitrary
frame of reference. Let $S$ and $T$ be new coordinates introduced as
follows:\begin{equation}S=-z\cosh\psi+t\sinh\psi;\enskip
T=t\cosh\psi-z\sinh\psi.\label{lorentz}\end{equation}Here the parameter
$\psi$ labels inertial reference frames moving along the $z$-axis.

Consider the following triplet of 2-forms:\begin{equation}
\pi^1=P(T)dT\wedge dZ-Q(T)\rho d\rho\wedge d\varphi\label{picyl}
\end{equation}$$\pi^2=p(T)dZ\wedge d\rho+q(T)\rho dT\wedge d\varphi$$
$$\pi^3=p(T)\rho d\varphi\wedge dZ+q(T)dT\wedge d\rho.$$
Similarly to the previous case the Hodge conjugation acts as the
substitution (\ref{hodge}). Inserting this into the first structure
equation (\ref{str1}) and source-free Yang-Mills equation (\ref{yme}) with
the connection 1-form $\alpha^a$:\begin{equation}
\alpha^1=f(T)dZ-d\varphi;\quad\alpha^2=g(T)\rho d\varphi;\quad
\alpha^3=g(T)d\rho.\label{alphacyl}\end{equation}gives the following:
\begin{equation}-Q'+2gq=0,\quad p'=Pg+FQ,\quad-P'+2gp=0,\quad-q'-gQ+pf=0.
\label{cylfour}\end{equation}The curvature 2-form corresponding to the
connection (\ref{alphacyl})is
$$K^1=f'dT\wedge dZ-g^2\rho d\rho\wedge d\varphi,$$
$$K^2=g'\rho dT\wedge d\varphi-fgdZ\wedge d\rho,$$
$$K^3=g'dT\wedge d\rho-fg\rho d\varphi\wedge dZ.$$Inserting this and
the expressions (\ref{picyl}) into the equation (\ref{fieldeq}) gives:
$$f'=m^2P,\quad g^2=m^2Q,\quad g'=m^2q,\quad -fg=m^2p.$$
Now we introduce the new function $h=f/\sqrt{2}$ and henceforth deal
only with the functions $g$ and $h$ which are\begin{equation}
g=m\sqrt{Q}\quad h=mp/2\sqrt{Q}.\label{gh}\end{equation}
The complete sextet of independent equations including (\ref{cylfour})
in this case is$$h''+2g^2h=0,\quad g''+g^3+2gh^2=0,\quad h'=m^2P\sqrt{2},$$
$$g^2=m^2Q,\quad g'=m^2Q,\quad gh=m^2p\sqrt{2}.$$Apparently, the first two
of them are to solved separately and all the rest define the explicit form
of the connection (\ref{alphacyl}). The independent equations
\begin{equation}h''+2g^2h=0,\quad g''+g^3+2gh^2=0\label{uniform}
\end{equation}draw a classical mechanical problem with Hamiltonian
\begin{equation}H=\frac{1}{2}(g'^2+h'^2)+\frac{g^4}{4}+\frac{g^2h^2}{2}.
\label{hamilt}\end{equation}In the case $h=0$ the system is well studied
and solutions in Jacobi elliptic functions with
$g=m\cdot sd(mT\enskip|\enskip1/2)$ are found long ago (Actor, 1979).   
Note that in the case $h=0$ there exists the only solution with mass
$m$, i.\,e., transformations (\ref{lorentz}) change the argument of
the elliptic function $mT$ into $\omega t-kz$ where
$\omega=m\cosh\psi$ and $k=m\sinh\psi$, thus, $\omega^2-k^2=m^2$. The
field equations have been reduced to the system (\ref{uniform}) and
the Hamiltonian (\ref{hamilt}) in our work (Turakulov 3)..
\section{The field of spherical wave}

The coordinate system used below is introduced in the framework of
the standard spherical coordinates $\{t,r,\theta,\varphi\}$ by
the following substitutions $\{\zeta,\eta,\theta,\varphi\}$:
$$t=\zeta\cdot\cosh\eta;\quad r=\zeta\cdot\sinh\eta;$$
$$\zeta=\sqrt{t^2-r^2};\quad\eta=arctanh(r/t).$$It is seen that
$\zeta$ is equal to interval between a given point and the cone point
$t=r=0$, and the surface $\zeta=0$ is the light cone. The coordinate
$\eta$ labels 3-dimensional cones with timelike generatrixes orthogonal
to the surfaces $\zeta=const$. The extremal values of $\eta$, namely,
$\eta=0$ and $\eta=\infty$ correspond to timelike straight line $r=0$,
which, meanwhile, could be chosen arbitrarily, and to the light cone
itself respectively.

The metric of this coordinate system is$$g=d\zeta\otimes d\zeta-
\zeta^2[d\eta\otimes d\eta+\sinh^2\eta(d\theta\otimes d\theta+\sin^2\theta
d\varphi\otimes d\varphi)]$$and as it is orthogonal the unit 4-form is
$[1-(\eta^2/4)]^{-3/2}\zeta^3\eta^2\sin\theta
d\zeta\wedge d\eta\wedge d\theta\wedge d\varphi$. Therefore,
non-zero components of Levi-Civita symbol are
$\varepsilon_{\zeta\eta\theta\varphi}=\zeta^3\sinh^2\eta\sin\theta$
Levi-Civita symbol components will be used when completing the
$\ast$-conjugating (Warner, 1983) of the 2-form of Yang-Mills field strength.

Consider the following triplet of 2-forms:\begin{equation}
\pi^1=p(\zeta)\zeta d\zeta\wedge d\eta-q(\zeta)\zeta^2\sinh^2\eta
\sin\theta d\theta\wedge d\varphi\label{pisw}
\end{equation}$$\pi^2=p(\zeta)\zeta\sinh\eta d\zeta\wedge d\theta
+q(\zeta)\zeta^2\sinh\eta\sin\theta d\varphi\wedge d\eta$$
$$\pi^3=p(\zeta)\zeta\sinh\eta\theta d\zeta\wedge d\varphi+
q(\zeta)\zeta^2\sinh\eta d\eta\wedge d\theta.$$ Here $P=p$ and $Q=q$
due to more wide symmetry group than that encountered in previous
cases. In this case the field is not only spherically symmetric but also
Lorentz-invariant. As usual, the Hodge conjugation acts as the
substitution $p\rightarrow q,\enskip\rightarrow q$. Now solving the first
structure equation (\ref{str1}) in combination with the source-free
Yang-Mills equation (\ref{yme}) with the connection 1-form $\alpha^a$:
$$\alpha^1=A(\zeta)d\eta+\cos\theta d\varphi;\quad
\alpha^2=A(\zeta)\sinh\eta d\theta-\cosh\eta\sin\theta d\varphi;$$ $$
\alpha^3=\cosh\eta d\theta+A(\zeta)\sinh\eta\sin\theta d\varphi$$
and the corresponding curvature form which is
$$K^1=A'd\zeta\wedge d\eta+
(1+A^2)\sinh^2\eta\sin\theta d\theta\wedge d\varphi$$
$$K^2=A'\sinh\eta d\zeta\wedge d\theta+
(1+A^2)\sinh\eta\sin\theta d\varphi\wedge d\eta$$
$$K^3=A'\sinh\eta\sin\theta d\zeta\wedge d\varphi+
(1+A^2)\sinh\eta d\eta\wedge d\theta$$
gives the following (we put $m=1$):$$A'=p\zeta,\quad(1+A^2)=q\zeta^2,
\quad(\zeta^2q)'=2Ap\zeta,\quad(\zeta^2p)'=-2Aq\zeta.$$Consequently,
the function $A$ satisfies the equation$$\zeta(\zeta A')'+2A(1+A^2)=0.$$
As seen from this equation the solutions are forminvariant under Lorentz
transformations as well as the equation itself. The substitutions
$s=\ln\zeta,\enskip A(\zeta)=y(s)$ reduce this equation to the form of
equations of motion of a point mass in classical mechanics, i.e., the
well-known Newton equation:$$y''+2y(1+y^2)=0.$$Now, applying the standard
methods we find that the function $y$ satisfies the following equation:
$$y'=\sqrt{c_1^2-(1+y^2)^2}$$whose solution is one of Jacobi elliptical
functions (Abramowicz and Stegun, 1964):$$y=sd\left((s+c_2)
\sqrt{\frac{2c_1}{c_1^2-1}}|\quad\frac{c_1-1}{2c_1}\right)$$
This result was obtained in our work (Turakulov 3, 1995).\section{Spin of
plane wave}

As a rule, when dealing with linear classical fields one determines spin
of the field as the number of vectorial indices the field has accounting
spinorial index as half of vectorial one. They are equal indeed if the
field is linear. In fact, spin of a classical field is defined as integral
of spin density over the space. In the case of plane waves the density is
proportional to number of indicies. The difference between linear and
non-linear fields arising in the case under consideration can be seen in
the following example. The usual solutions of Maxwell equations that
describe linearly polarized plane waves correspond zero spin density
because they represent a mixed state with $+1$ and $-1$ helicities. To
have the proper value with spin $1$ To have the proper value with spin
$1$ one composes a complex linear combination of two plane waves which
has circular polarization, and considers it as a pure state with certain
helicity. But this is possible only for linear fields. If the field is
non-linear no such linear combinations satisfy the original equation.

In this section we evaluate the spin density of plane wave. Applying
the N\"other thorem gives the spin density in the form$$S^t_{ij}=
\frac{\partial L}{\partial D_tA^a_l}\delta_k{}^l{}_{ij}A^a_l=
\delta_k{}^l{}_{ij}\pi_a^{tk}A^a_l$$where we used the denotion
$\delta_k{}^l{}_{ij}$ for $(1/2)(\eta_{ki}\delta^l_j-
\eta_{kj}\delta^l_i$. Thus,$$S^t_{ij}=(1/2)[(\pi_a)^t{.j}A^a_i-
(\pi_a)^t_{.j}A^a_i]$$Now, a straightforward evaluation for the plane
wave solution (\ref{picyl}) and (\ref{alphacyl}) shows that the only
non-zero component of the spin density is $S^t_{.z\varphi}=-f'$.
Apparently, this density gives zero total spin because locally it is
everywhere pointed out from the axis. Consequently, plane waves of the
field are spinless.\section{Acknowledgments}

The author thanks Prof. H.-D.\,Doebner for hospitality in XXI
International Colloquium on Group Theoretic Methods in Physics and
Prof. B.\,F.\,Schutz for hospitality in Albert Einstein Institute,
Potsdam and for interesting discussions.

{\Large\bf References}\\
M Abramowitz and I Stegun, (1964) {\it Handbook of Mathematical
Functions} (National Bureau of Standards)\\
A Actor (1979) {\it Rev Mod Phys} {\bf 51} 461\\
M Atiah and N Hitchin, {\it The Geometry and Dynamics
of Magnetic Monopoles} (1988) (Princeton University Press, Princeton) \\
M Calvo (1977){\it Phys Rev} {\bf D15} 1733\\
L Castillejo, M Kugler  (1980) {\it Phys Lett} {\bf B93} 61\\
T Eguchi, P B Gilkey and A J Hanson, (1980) {\it Phys. Reports}
{\bf C66} 214\\
M B Halpern {\it Phys Rev} (1977) {\bf D16} 1798\\
J P Hsu, E Mac (1977) {\it J. Math. Phys.}{\bf  18}  100\\
R Roskies (1977) {\it Phys Rev} {\bf D15} 1731\\
H Ryder (1985){\it Quantum Field Theory} (Cambridge University
Press, Cambridge)\\
M Socolowski (1991) {\it J Math Phys} {\bf 32}  2522\\
A Trautman (1979) {\it Czech J Phys} {\bf B29} 107\\
Z Y Turakulov 1 (1995) {\it Tr J of Physics} {\bf 19} 632\\ 
Z Y Turakulov 2 (1995) {\it Tr J of Physics} {\bf 18} 1346\\
Z Y Turakulov 3 (1995) {\it Tr J of Physics} {\bf 19} 635\\
F.\,Warner, (1983) {\it Foundations of Differentiable Manifolds and
Lie Groups} (Springer, New York)\\
D Wilkinson, A Goldhaber (1977) {\it Phys Rev} {\bf D16} 1221\\
A C T Wu, T T Wu (1974) {\it J Math Phys} {\bf 15} 53\\
T T Wu, C N Yang (1975) {\it Phys Rev} {\bf D12} 3843\\
P Yasskin (1980) in: {\it Gauge Theories} Eds.\, J P Harnad
and S Shnider (Springer, Berlin) 154\\
\end{document}